# Linearly edge-reinforced random walks

Franz Merkl[1] and Silke W. W. Rolles[2]

*University of Munich and University of Bielefeld*

**Abstract:** We review results on linearly edge-reinforced random walks. On finite graphs, the process has the same distribution as a mixture of reversible Markov chains. This has applications in Bayesian statistics and it has been used in studying the random walk on infinite graphs. On trees, one has a representation as a random walk in an independent random environment. We review recent results for the random walk on ladders: recurrence, a representation as a random walk in a random environment, and estimates for the position of the random walker.

## 1. Introduction

Consider a locally finite graph $G = (V, E)$ with undirected edges. All the edges are assigned weights which change in time. At each discrete time step, the *edge-reinforced random walker* jumps to a neighboring vertex with probability proportional to the weight of the traversed edge. As soon as an edge is traversed, its weight is increased by one.

The model was introduced by Diaconis in 1986 (see [11] and [3]). Mike Keane made the model popular in The Netherlands. In particular, Mike Keane introduced Silke Rolles to the intriguing questions raised by this model. One fundamental question, asked by Diaconis, concerns *recurrence*: Do almost all paths of the edge-reinforced random walk visit all vertices infinitely often? By a Borel-Cantelli argument, this is equivalent to the following question: Does the edge-reinforced random walker return to the starting point infinitely often with probability one? For all dimensions $d \geq 2$, it is an open problem to prove or disprove recurrence on $\mathbb{Z}^d$. Only recently, recurrence of the edge-reinforced random walk on ladders has been proven.

The present article focuses on *linear edge-reinforcement* as described above. In the past two decades, many different reinforcement schemes have been studied. We briefly mention some of them: *Vertex-reinforced random walk* was introduced by Pemantle [25] and has been further analyzed by Benaïm [1], Dai ([5], [6]), Pemantle and Volkov ([27], [38]), and Tarrès [35]. Sellke [30] and Vervoort [37] proved recurrence results for *once-reinforced random walks* on ladders; Durrett, Kesten, and Limic [13] analyzed the process on regular trees. Recurrence questions for *reinforced random walks of sequence type* were studied by Davis ([7], [8]) and Sellke [31]. Takeshima ([33], [34]) studied hitting times and recurrence for *reinforced random walk of matrix type*. Limic [18] proved a localization theorem for superlinear reinforcement. *Weakly reinforced random walks* in one dimension were studied by

[1]Mathematical Institute, University of Munich, Theresienstr. 39, D-80333 Munich, Germany, e-mail: merkl@mathematik.uni-muenchen.de
[2]Department of Mathematics, University of Bielefeld, Postfach 100131, D-33501 Bielefeld, Germany, e-mail: srolles@math.uni-bielefeld.de

*AMS 2000 subject classifications:* primary 82B41; secondary 60K35, 60K37.

*Keywords and phrases:* reinforced random walk, random environment, recurrence, Bayesian statistics.

Tóth [36]. *Directionally reinforced random walks* are the subject of Mauldin, Monticino, and von Weizsäcker [20] and Horváth and Shao [14]. Othmer and Stevens [23] suggested linearly edge-reinforced random walks as a simple model for the gliding of myxobacteria. These bacteria produce a slime and prefer to move on the slime trail produced earlier. Overviews with a different emphasis than the present article have been written by Davis [9] and Pemantle [26].

The article is organized as follows: In Section 2, we give a formal definition of the edge-reinforced random walk. On finite graphs, the edge-reinforced random walk has the same distribution as a random walk in a dependent random environment. This representation together with related limit theorems is presented in Section 3. The edge-reinforced random walk on finite graphs gives rise to a family of prior distributions for reversible Markov chains. This application to statistics is the content of Section 4. A characterization of the process (which is also applied in the Bayesian context) is given in Section 5. In Section 6, we state results for the process on acyclic graphs. Recently, progress has been made in understanding the edge-reinforced random walk on graphs of the form $\mathbb{Z} \times \{1, \ldots, d\}$ and, more generally, on $\mathbb{Z} \times T$ when $T$ is a finite tree. Section 7 reviews these recent results and presents some simulations.

## 2. Formal description of the model

Let $G = (V, E)$ be a locally finite undirected graph. Every edge is assumed to have two different endpoints; thus there are no direct loops. We identify an edge with the set of its endpoints. Formally, the edge-reinforced random walk on $G$ is defined as follows: Let $X_t$ denote the random location of the random walker at time $t$. Let $a_e > 0$, $e \in E$. For $t \in \mathbb{N}_0$, we define $w_t(e)$, the weight of edge $e$ at time $t$, recursively as follows:

$$w_0(e) := a_e \quad \text{for all } e \in E, \tag{2.1}$$

$$w_{t+1}(e) := \begin{cases} w_t(e) + 1 & \text{for } e = \{X_t, X_{t+1}\} \in E, \\ w_t(e) & \text{for } e \in E \setminus \{\{X_t, X_{t+1}\}\}. \end{cases} \tag{2.2}$$

Let $P_{v_0,a}$ denote the distribution of the edge-reinforced random walk on $G$ starting in $v_0$ with initial edge weights equal to $a = (a_e)_{e \in E}$. The distribution $P_{v_0,a}$ is a probability measure on $V^{\mathbb{N}_0}$, specified by the following requirements:

$$X_0 = v_0 \quad P_{v_0,a}\text{-a.s.}, \tag{2.3}$$

$$P_{v_0,a}(X_{t+1} = v | X_i, i = 0, 1, \ldots, t) = \begin{cases} \dfrac{w_t(\{X_t, v\})}{\sum_{\{e \in E: X_t \in e\}} w_t(e)} & \text{if } \{X_t, v\} \in E, \\ 0 & \text{otherwise.} \end{cases} \tag{2.4}$$

## 3. Reinforced random walk on finite graphs

Throughout this section, we assume the graph $G$ to be finite.

For $t \in \mathbb{N}$ and $e \in E$, set

$$k_t(e) := w_t(e) - a_e \quad \text{and} \quad \alpha_t(e) := \frac{k_t(e)}{t}. \tag{3.1}$$

In particular, $\alpha_t(e)$ denotes the proportion of crossings of the edge $e$ up to time $t$. The random vector $\alpha_t := (\alpha_t(e))_{e \in E}$ takes values in the simplex

$$\Delta := \left\{ (x_e)_{e \in E} \in (0, \infty)^E : \sum_{e \in E} x_e = 1 \right\}. \tag{3.2}$$

For $x \in (0, \infty)^E$ and $v \in V$, we define

$$x_v := \sum_{e \in E : v \in e} x_e. \tag{3.3}$$

Now view the graph $G$ as an electric network, and consider the space $H_1$ of all electric current distributions $(y_e)_{e \in E}$ on the graph, such that Kirchhoff's vertex rule holds at all vertices: For every vertex, the sum of all ingoing currents should equal the sum of all outgoing currents. More formally, we proceed as follows: We assign a counting direction to every edge, encoded in a signed incidence matrix $s = (s_{ve})_{v \in V, e \in E} \in \{-1, 0, 1\}^{V \times E}$. Here, $s_{ve} = +1$ means that $e$ is an ingoing edge into the vertex $v$; $s_{ve} = -1$ means that $e$ is an outgoing edge from vertex $v$, and $s_{ve} = 0$ means that the vertex $v$ is not incident to the edge $e$. Then, for a given current distribution $y = (y_e)_{e \in E} \in \mathbb{R}^E$, the current balance at any vertex $v \in V$ equals

$$(sy)_v = \sum_{e \in E} s_{ve} y_e, \tag{3.4}$$

and the space of all current distributions satisfying Kirchhoff's vertex rule equals

$$H_1 = \text{kernel}(s) = \{ y \in \mathbb{R}^E : sy = 0 \}$$
$$= \left\{ (y_e)_{e \in E} \in \mathbb{R}^E : \sum_{e \in E} s_{ve} y_e = 0 \text{ for all } v \in V \right\}. \tag{3.5}$$

In other words, $H_1$ is the first homology space of the graph $G$. Note that $\dim H_1 = |E| - |V| + 1$ by Euler's rule.

Interpreting weights $x = (x_e)_{e \in E} \in (0, \infty)^E$ as electric conductivities of the edges, the power consumption of the network required to support a current distribution $y \in H_1$ equals

$$A_x(y) = \sum_{e \in E} \frac{y_e^2}{x_e}.$$

(Of course, in order to drive this current in a real network, appropriate batteries must be wired into the edges. We ignore this fact.) Note that $A_x$ is a quadratic form on $H_1$. We will also need its determinant $\det A_x$, which is defined as the determinant of the matrix representing $A_x$ with respect to any $\mathbb{Z}$-basis of the lattice $H_1 \cap \mathbb{Z}^E$. Note that the determinant $\det A_x$ does not depend on the choice of the basis, since any base change has determinant $\pm 1$. We endow the homology space $H_1$ with the Lebesgue measure, normalized such that the unit cell spanned by any $\mathbb{Z}$-basis of $H_1 \cap \mathbb{Z}^E$ gets volume 1.

It turns out that the sequence $(\alpha_t)_{t \in \mathbb{N}}$ converges almost surely to a random limit. Surprisingly, the limiting distribution can be determined explicitly:

**Theorem 3.1.** *The sequence $(\alpha_t)_{t \in \mathbb{N}}$ converges $P_{v_0, a}$-almost surely. The distribution of the limit is absolutely continuous with respect to the surface measure on $\Delta$*

*with density given by*

$$\phi_{v_0,a}(x) = Z_{v_0,a}^{-1} \frac{\prod_{e \in E} x_e^{a_e - \frac{1}{2}}}{x_{v_0}^{\frac{a_{v_0}}{2}} \prod_{v \in V \setminus \{v_0\}} x_v^{\frac{a_v+1}{2}}} \sqrt{\det A_x}, \qquad x = (x_e)_{e \in E} \in \Delta, \qquad (3.6)$$

*where $Z_{v_0,a}$ denotes a normalizing constant.*

Theorem 3.1 was discovered by Coppersmith and Diaconis [3]. A proof of the result was independently found by Keane and Rolles ([16], Theorem 1). In the special case of a triangle, a derivation of Theorem 3.1 was published by Keane [15]. The normalizing constant $Z_{v_0,a}$ is explicitly known; see [16].

Let $\mathcal{S}$ denote the set of all subtrees $T$ in $G$, viewed as set of edges $T \subseteq E$, which need not necessarily visit all vertices. The set of all *spanning* trees in $G$ is denoted by $\mathcal{T}$; it is a subset of $\mathcal{S}$. Using a matrix tree theorem (see e.g. [19], page 145, theorem 3'), one can rewrite the determinant of $A_x$ as a sum over all spanning trees. This yields the following representation of the density:

$$\phi_{v_0,a}(x) = Z_{v_0,a}^{-1} \frac{\prod_{e \in E} x_e^{a_e - 1}}{x_{v_0}^{\frac{a_{v_0}}{2}} \prod_{v \in V \setminus \{v_0\}} x_v^{\frac{a_v+1}{2}}} \sqrt{\sum_{T \in \mathcal{T}} \prod_{e \in T} x_e}. \qquad (3.7)$$

Given a path $\pi = (v_0, v_1, \ldots, v_t)$ of vertices in $G$, we define its corresponding "chain" $[\pi] = ([\pi]_e)_{e \in E} \in \mathbb{R}^E$ as follows: Imagine an electric current of size 1 entering the network $G$ at the starting vertex $v_0$, flowing along $\pi$, and then leaving the network at the last vertex $v_t$. Let $[\pi] \in \mathbb{R}^E$ denote the corresponding current distribution. More formally, for $e \in E$, $[\pi]_e$ equals the number of times the path $\pi$ traverses the edge $e$ in counting direction minus the number of times the path $\pi$ traverses the edge $e$ in opposite direction to the counting direction. Note that $[\pi] \notin H_1$ unless the path $\pi$ is closed, since Kirchhoff's vertex rule is violated at the starting point $v_0$ and at the last vertex $v_t$. Now, for any time $t$, let

$$\pi_t = (X_0, X_1, \ldots, X_t)$$

denote the random path the reinforced random walker follows up to time $t$. We define

$$\beta_t = \frac{1}{\sqrt{t}}[\pi_t] \in \mathbb{R}^E \qquad (3.8)$$

to be the corresponding rescaled (random) current distribution. Indeed, the diffusive scale $\sqrt{t}$ turns out to be appropriate for studying the random currents $[\pi_t]$ in the limit as $t \to \infty$. Note that in general $\beta_t \notin H_1$, due to the violation of Kirchhoff's rule at the starting vertex $X_0$ and the last vertex $X_t$. However, any weak limit of $\beta_t$ as $t \to \infty$, if it exists, must be supported on the homology space $H_1 \subseteq \mathbb{R}^E$. Indeed, the violation of Kirchhoff's rule at any vertex is not larger than $t^{-1/2}$, which is negligible as $t \to \infty$.

For a (random) path $\pi_t = (X_0, \ldots, X_t)$ and any vertex $v$ among $\{X_0, \ldots, X_{t-1}\}$, let $e^{\text{last exit}}(v, \pi_t)$ denote the edge by which $v$ is left when it is visited by $\pi_t$ for the last time. Let $T_t^{\text{last exit}} \in \mathcal{S}$ denote the random tree graph consisting of all the edges $e^{\text{last exit}}(v, \pi_t)$, $v \in \{X_0, \ldots, X_{t-1}\}$. Since the reinforced random walk on the finite graph $G$ visits every vertex almost surely, $T_t^{\text{last exit}} \in \mathcal{T}$ holds for all large $t$ almost surely.

**Theorem 3.2.** *The sequence $(\alpha_t, \beta_t, T_t^{\text{last exit}})_{t\in\mathbb{N}}$ converges weakly in $\mathbb{R}^E \times \mathbb{R}^E \times \mathcal{S}$. The limiting distribution $\mathbb{Q}$ is supported on the subset $\Delta \times H_1 \times \mathcal{T}$ of $\mathbb{R}^E \times \mathbb{R}^E \times \mathcal{S}$. It is absolutely continuous with respect to the product of the surface measure on $\Delta$, the Lebesgue measure on $H_1$, and the counting measure on $\mathcal{T}$. Its density is given by*

$$\Phi_{v_0,a}(x,y,T) = \tilde{Z}_{v_0,a}^{-1} \frac{\prod_{e\in E} x_e^{a_e - \frac{3}{2}}}{x_{v_0}^{\frac{a_{v_0}}{2}} \prod_{v\in V\setminus\{v_0\}} x_v^{\frac{a_v+1}{2}}} \left(\prod_{e\in T} x_e\right) \exp\left(-\frac{1}{2} A_x(y)\right), \quad (3.9)$$

*where $\tilde{Z}_{v_0,a}$ denotes a normalizing constant.*

Theorem 1 of [16] states weak convergence of $(\alpha_t, \beta_t)_{t\in\mathbb{N}}$ to the limiting distribution with density

$$\sum_{T\in\mathcal{T}} \Phi_{v_0,a}(x,y,T). \quad (3.10)$$

However, the proof implicitly contains the proof of Theorem 3.2.

Recall that the transition probabilities of any irreducible reversible Markov chain on $G$ can be described by weights $x = (x_e)_{e\in E}$, $x_e \geq 0$, on the edges of the graph; the probability to traverse an edge is proportional to its weight. More precisely, denoting the distribution of the Markov chain induced by the edge weights $x$ with starting vertex $v_0$ by $Q_{v_0,x}$, one has

$$Q_{v_0,x}(X_{t+1} = v' | X_t = v) = \frac{x_{\{v,v'\}}}{x_v}, \quad (3.11)$$

whenever $\{v,v'\}$ is an edge; the weight $x_v$ is defined in (3.3).

The following representation of the edge-reinforced random walk on a finite graph as a mixture of reversible Markov chains is shown in Theorem 3.1 of [28].

**Theorem 3.3.** *For any event $B \subseteq V^{\mathbb{N}_0}$, one has*

$$P_{v_0,a}((X_t)_{t\in\mathbb{N}_0} \in B) = \int_{\Delta \times H_1 \times \mathcal{T}} Q_{v_0,x}(B) \, \mathbb{Q}(dx\,dy\,dT), \quad (3.12)$$

*where $\mathbb{Q}$ denotes the limiting measure from Theorem 3.2. In other words, the edge-reinforced random walk on any finite graph $G$ has the same distribution as a random walk in a random environment. The latter is given by random weights on the edges, distributed according to the limiting distribution of $(\alpha_t)_{t\in\mathbb{N}}$, namely the distribution with density $\phi_{v_0,a}$ with respect to the Lebesgue measure on the simplex.*

This representation as a random walk in a random environment has been extremly useful in studying the edge-reinforced random walk on infinite ladders.

The proof of Theorem 3.3 relies on the fact that the edge-reinforced random walk is *partially exchangeable*: the probability that the process traverses a particular finite path $\pi = (v_0, \ldots, v_t)$ depends only on the starting point of $\pi$ and the number of transition counts of the undirected edges. If one knows that the process returns to the starting point infinitely often with probability one (which is the case for the edge-reinforced random walk on a finite graph), one can apply a de Finetti Theorem for Markov chains of Diaconis and Freedman [10]. A refinement for reversible Markov chains ([28], Theorem 1.1) yields a mixture of *reversible chains*.

## 4. An application to Bayesian statistics

Consider the following statistical situation: We observe $X_0 = v_0$, $X_1 = v_1$, ..., $X_t = v_t$ generated by a reversible Markov chain. The transition kernel $k(\cdot, \cdot)$ and the stationary measure $\nu$ are *unknown*. Let $V$ be the set of possible observations. We assume $V$ to be a *known* finite set. Furthermore, we assume that $k(v, v') > 0$ if and only if $k(v', v) > 0$. Hence, $V$ together with the set $E := \{\{v, v'\} : k(v, v') > 0\}$ defines a finite undirected graph. This graph is assumed to be *known*.

It is a natural question how to model this in the framework of Bayesian statistics. One is interested in "natural" prior distributions on the set of reversible Markov chains.

Because of reversibility, $\nu(v)k(v, v') = \nu(v')k(v', v)$ for all $v, v' \in V$. The distribution of a reversible Markov chain with transition kernel $k$, stationary distribution $\nu$, and starting point $v_0$ is given by $Q_{v_0,x}$ as defined in (3.11) with edge weights $x_{\{v,v'\}} = \nu(v)k(v, v')$. Thus, one can describe the prior distributions as measures on the set of possible edge weights, namely as measures on $\Delta$.

The following theorem was proved by Diaconis and Rolles (Proposition 4.1 of [12]):

**Theorem 4.1.** *Let $\sigma$ denote the Lebesgue measure on $\Delta$ and set $\mathbb{P}_{v_0,a} := \phi_{v_0,a} \, d\sigma$. The family*

$$\{\mathbb{P}_{v_0,a} \; : \; v_0 \in V, a = (a_e)_{e \in E} \in (0, \infty)^E\} \tag{4.1}$$

*of prior distributions is closed under sampling. More precisely, under the prior distribution $\mathbb{P}_{v_0,a}$ with observations $X_0 = v_0$, $X_1 = v_1$, ..., $X_t = v_t$, the posterior is given by $\mathbb{P}_{v_t,(a_e+k_t(e))_{e \in E}}$, where $k_t(e)$ denotes the number of $i$ with $\{v_i, v_{i+1}\} = e$.*

These prior distributions were further analyzed in [12]: They can be generalized to finite graphs with direct loops; thus one can include the case $k(v, v) > 0$ for some $v$. The set of linear combinations of the priors $\mathbb{P}_{v_0,a}$ is weak-star dense in the set of all prior distributions on reversible Markov chains on $G$. Furthermore, it is shown that these prior distributions allow to perform tests: several hypotheses are tested for a data set of length 3370 arising from the DNA sequence of the humane HLA-B gene, for instance $H_0$ : i.i.d.(unknown) versus $H_1$ : reversible Markov chain and $H_0$ : reversible Markov chain versus $H_1$ : full Markov. The tests are based on the Bayes factor $P(\text{data}|H_0)/P(\text{data}|H_1)$ which can be easily computed.

The priors $\mathbb{P}_{v_0,a}$ generalize the well-known Dirichlet priors. The latter are obtained as a special case for star-shaped graphs. The Dirichlet prior was characterized by W.E. Johnson; see [39]. In [12], a similar characterization is given for the priors $\mathbb{P}_{v_0,a}$; in this sense they are "natural".

## 5. A characterization of reinforced random walk

In this section, we review a characterization of the edge-reinforced random walk from [28]. We need some assumptions on the underlying graph $G$:

**Assumption 5.1.** *For all $v \in V$, degree$(v) \neq 2$. Furthermore, the graph $G$ is 2-edge-connected, i.e. removing an edge does not make $G$ disconnected.*

Let $P$ be the distribution of a nearest-neighbor random walk on $G$ such that the following hold:

**Assumption 5.2.** *There exists $v_0 \in V$ with $P(X_0 = v_0) = 1$.*

**Assumption 5.3.** For any admissible path $\pi$ of length $t \geq 1$ starting at $v_0$, we have $P((X_0, \ldots, X_t) = \pi) > 0$.

**Assumption 5.4.** The process $(X_t)_{t \in \mathbb{N}_0}$ with distribution $P$ is partially exchangeable.

For $t \in \mathbb{N}_0$, $v \in V$, and $e \in E$, we define

$$k_t(v) := |\{i \in \{0, \ldots, t\} : X_i = v\}|, \quad k_t(e) := |\{i \in \{1, \ldots, t\} : \{X_{i-1}, X_i\} = e\}|.$$

**Assumption 5.5.** For all $v \in V$ and $e \in E$, there exists a function $f_{v,e}$ taking values in $[0, 1]$ such that for all $t \in \mathbb{N}_0$

$$P(X_{t+1} = v | X_0, \ldots, X_t) = f_{X_t, \{X_t, v\}}(k_t(X_t), k_t(\{X_t, v\})).$$

In other words, the conditional distribution for the next step, given the past up to time $t$, depends only on the position $X_t$ at time $t$, the edge $\{X_t, v\}$ to be traversed, the local time accumulated at $X_t$, and the local time on the edge $\{X_t, v\}$.

It is not hard to see that an edge-reinforced random walk and a non-reinforced random walk starting at $v_0$ satisfy Assumptions 5.2–5.5. The following theorem is the content of Theorem 1.2 of [28]:

**Theorem 5.1.** *Suppose the graph $G$ satisfies Assumption 5.1. If $P$ is the distribution of a nearest-neighbor random walk on $G$ satisfying Assumptions 5.2–5.5, then, for all $t$,*

$$P(X_{t+1} = v | X_0, \ldots, X_t)$$

*agrees on the set $\{k_t(X_t) \geq 3\}$ with the corresponding conditional probability for an edge-reinforced random walk or a non-reinforced random walk starting at $v_0$.*

In this sense, the above assumptions characterize the edge-reinforced random walk. Theorem 5.1 is used to give a charactarization of the priors $\mathbb{P}_{v_0, a}$.

## 6. Reinforced random walk on acyclic graphs

The edge-reinforced random walk on *acyclic* graphs is much easier to analyze than on graphs with cycles. Let us briefly explain why: Consider edge-reinforced random walk on a tree. If the random walker leaves vertex $v$ via the neighboring vertex $v'$, then in case the random walker ever returns to $v$, the next time she does so, she *has to enter via the same edge* $\{v, v'\}$. Hence, if the random walker leaves $v$ via the edge $\{v, v'\}$, the weight of this edge will have increased by precisely 2, the next time the random walker arrives at $v$. Obviously, this is only true on an acyclic graph.

Due to this observation, instead of recording the edge weights, one can place *Polya urns* at the vertices of the tree. Each time the random walker is at location $v$, a ball is drawn from the urn $U(v)$ attached to $v$. Then, the ball is put back together with two balls of the same color. The different colors represent the edges incident to $v$; the number of balls in the urn equal the weights of the edges incident to $v$, observed at times when the random walker is at $v$. The initial composition of the urns is determined by the starting point and the initial edge weights of the reinforced random walk. The sequence of drawings from the urn $U(v)$ at $v$ is *independent* from the sequences of drawings from the urns $U(v')$, $v \neq v'$, at all other locations. Using de Finetti's theorem, one finds that the edge-reinforced random walk has the same distribution as a random walk in a random environment, where the latter is given by *independent* Dirichlet-distributed transition probabilities. This representation as a random walk in a random environment was observed by Pemantle [24]. It seems to be the most powerful tool to analyze the process on acyclic graphs. Pemantle used it to prove a phase transition in the recurrence/transience behavior on a binary tree. Later, Takeshima [33] characterized recurrence in one dimension for space-inhomogeneous initial weights. (For all initial weights being equal, recurrence in one dimension follows for instance from a well-known recurrence criterion for random walk in random environment; see e.g. [32].) Collevecchio [4] proved a law of large numbers and a central limit theorem for the edge-reinforced random walk on $b$-ary trees, under the assumption that $b$ is large.

For linear edge-reinforcement on arbitrary *directed* graphs, so called *directed-edge-reinforced random walk*, a similar correspondance with a random walk in an independent random environment can be established. Using this correspondance, recurrence for $\mathbb{Z} \times G$ for any finite graph $G$ was proved by Keane and Rolles in [17]. A different criterion for recurrence and transience of a random walk in an independent random environment on a strip was established by Bolthausen and Goldsheid [2].

## 7. Reinforced random walk on ladders

Studying reinforced random walks on graphs of the form $\mathbb{Z} \times \{1, \ldots, d\}$ seems to be a challenging task. Only recently, recurrence results were obtained. The following result was proved in [22].

**Theorem 7.1.** *For all $a > 3/4$, the edge-reinforced random walk on $\mathbb{Z} \times \{1, 2\}$ with all initial weights equal to $a$ is recurrent.*

The result was extended in [29], under the assumption the initial weight $a$ is large:

**Theorem 7.2.** *Let $G$ be a finite tree. For all large enough $a$, the edge-reinforced random walk on $\mathbb{Z} \times G$ with all initial weights equal to $a$ is recurrent.*

In particular, this applies to ladders $\mathbb{Z} \times \{1, 2, \ldots, d\}$ of any finite width $d$.

In the following, we consider the process on $\mathbb{N}_0 \times G$ or $\mathbb{Z} \times G$ with a finite tree $G$. We always assume that all initial edge weights are equal to the same large enough constant $a$ and the random walk starts in a vertex $\mathbf{0}$ at level 0.

Just as edge-reinforced random walk on *finite* graphs, the edge-reinforced random walk on *infinite* ladders turns out to be equivalent to a random walk in a random environment, as studied in Section 3: In particular, there is an infinite-volume analogue to Theorem 3.3 for the *infinite* graphs $\mathbb{N}_0 \times G$ and $\mathbb{Z} \times G$. Here, the law $\mathbb{Q}$ of the random environment in Theorem 3.3 gets replaced by an infinite-volume Gibbs measure, which we also denote by $\mathbb{Q}$. Also in the infinite-volume setup, the fractions $\alpha_t \in \mathbb{R}^E$ of times spent on the edges converge almost surely to random weights $x \in \mathbb{R}^E$ as $t \to \infty$, just as in Theorem 3.1. These random weights $x$ are governed by the infinite-volume Gibbs measure $\mathbb{Q}$. It turns out that $\mathbb{Q}$-almost surely, the random weights decrease exponentially in space, i.e.

$$\mathbb{Q}\text{-a.s.} \qquad \limsup_{|e| \to \infty} \frac{1}{|e|} \log x_e \leq -c(a, G) \tag{7.1}$$

with a deterministic constant $c(G, a) > 0$. Here, $|e|$ denotes the distance of an edge $e$ from the starting point. For more details, see [21].

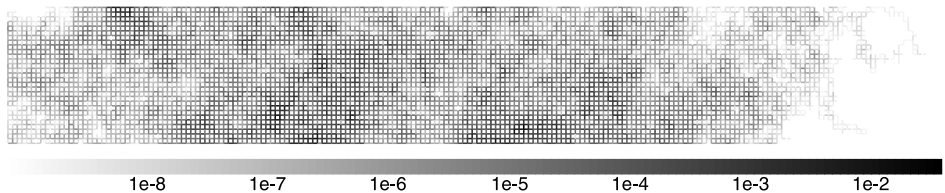

FIG 1. *Fraction of time spent on the edges on a logarithmic scale. Time horizon $= 10^9$, initial weight $a = 1$.*

This exponential decay of the weights $x$ can be also observed in simulations: For the ladder $\mathbb{Z} \times \{1, 2, \ldots, 30\}$, the simulation in Figure 1 shows one typical sample for the fraction of time $\alpha_t$ spent on any edge up to time $t = 10^9$. The fractions $\alpha_t$ are displayed logarithmically as gray scales. The starting point is located in the center of the picture. The time $t = 10^9$ in the simulation is already so large that this $\alpha_t$ is a good approximation for a typical sample of the random weights $x$.

Note that the fractions $\alpha_t$ in the simulation vary over many orders of magnitude, and that they decrease roughly exponentially as one gets farther away from the starting point.

Entropy estimates and deformation arguments from statistical mechanics are used in [21] to derive the exponential decay of the weights.

As a consequence, one obtains the following estimates for the position of the random walker, also proved in [21]:

**Theorem 7.3.** *There exist constants $c_1, c_2 > 0$, depending only on $G$ and $a$, such that for all $t, n \in \mathbb{N}_0$, the following bound holds:*

$$P_{\mathbf{0},a}(|X_t| \geq n) \leq c_1 e^{-c_2 n}. \tag{7.2}$$

Note that the bounds are uniform in the time $t$. This is different from the behavior of simple random walk, which has fluctuations of order $\sqrt{t}$.

**Corollary 7.1.** *There exists a constant $c_3 = c_3(G, a) > 0$ such that $P_{\mathbf{0},a}$-a.s.,*

$$\max_{s=0,\ldots,t} |X_s| \leq c_3 \ln t \qquad \text{for all } t \text{ large enough.} \tag{7.3}$$

A simulation shown in Figure 2 illustrates this corollary: For one sample path of an edge-reinforced random walk on the ladder $\mathbb{Z} \times \{1, 2, \ldots, 30\}$ with initial weight $a = 1$, the farthest level reached so far, $\max_{s=0,\ldots,t} |X_s|$, is displayed as a function of $t$. Note that the time $t$ is plotted on a logarithmic scale.

Figure 3 shows a simulation of reinforced random walk on $\mathbb{Z}^2$ with initial weights $a = 1$. It is still unknown whether this reinforced random walk is recurrent. Simulations show that there are random regions in $\mathbb{Z}^2$, maybe in far distance from the starting point, which are visited much more frequently than other regions closer to the starting point. It remains unclear whether more and more extreme "favorable regions" arise farther away from the origin. Thus, the recurrence problem in $\mathbb{Z}^2$ remains open, even on a heuristic level.

### Acknowledgment

We would like to thank Frank den Hollander and Guido Elsner for carefully reading our manuscript.

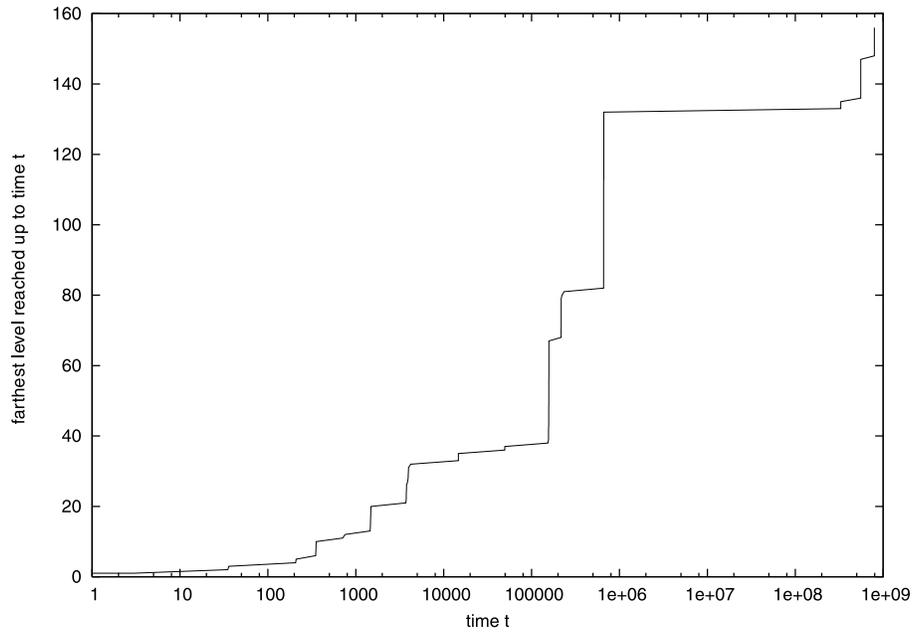

FIG 2. *Maximal distance from the starting point as a function of time. Initial weight $a = 1$.*

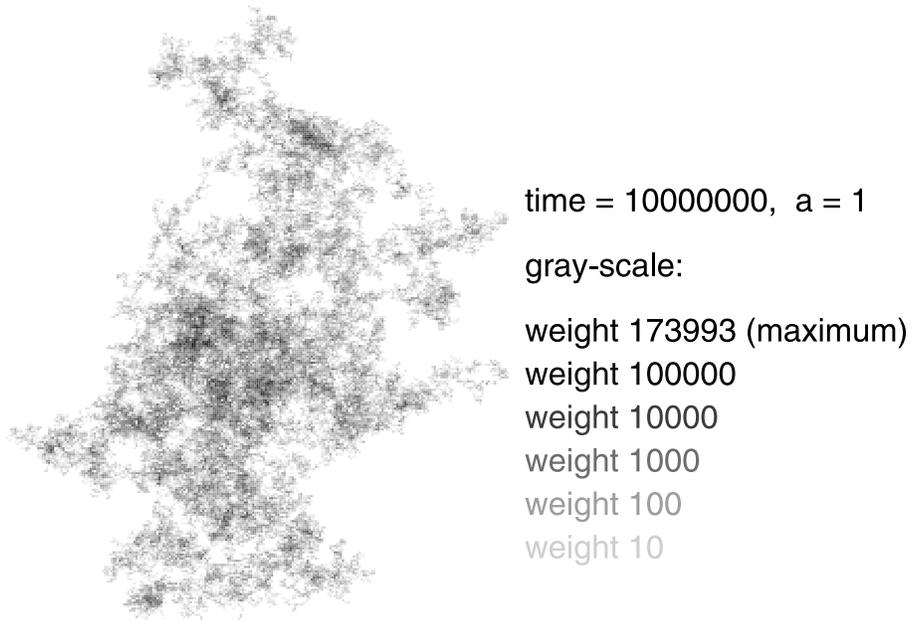

FIG 3. *Reinforced random walk on $\mathbb{Z}^2$.*